# SADDLEPOINT APPROXIMATION FOR MOMENT GENERATING FUNCTIONS OF TRUNCATED RANDOM VARIABLES[1]

By Ronald W. Butler and Andrew T. A. Wood

*Colorado State University and University of Nottingham*

We consider the problem of approximating the moment generating function (MGF) of a truncated random variable in terms of the MGF of the underlying (i.e., untruncated) random variable. The purpose of approximating the MGF is to enable the application of saddlepoint approximations to certain distributions determined by truncated random variables. Two important statistical applications are the following: the approximation of certain multivariate cumulative distribution functions; and the approximation of passage time distributions in ion channel models which incorporate time interval omission. We derive two types of representation for the MGF of a truncated random variable. One of these representations is obtained by exponential tilting. The second type of representation, which has two versions, is referred to as an exponential convolution representation. Each representation motivates a different approximation. It turns out that each of the three approximations is extremely accurate in those cases "to which it is suited." Moreover, there is a simple rule of thumb for deciding which approximation to use in a given case, and if this rule is followed, then our numerical and theoretical results indicate that the resulting approximation will be extremely accurate.

**1. Introduction.**

1.1. *Saddlepoint methods.* Saddlepoint methods provide approximations to densities and probabilities which are very accurate in a wide variety of settings. This accuracy is seen not only in numerical work, but also in theoretical calculations. In particular, it is often the case that relative errors of these approximations stay bounded in the extreme tails, a desirable property which is not shared by most other types of approximation used in statistics.

Received October 2002; revised February 2004.
[1]Supported in part by EPSRC.

For development and discussion of saddlepoint methodology and related methods, see Daniels (1954) for details of the density approximation; Barndorff-Nielsen and Cox (1989, 1994) for applications to inference; Lugannani and Rice (1980), Temme (1982) and Daniels (1987) for discussion of a tail area approximation which has uniform relative error, and Skovgaard (1987) for a conditional version of this approximation; and Reid (1988) for a review of saddlepoint techniques.

Saddlepoint approximations are constructed by performing various operations on the moment generating function (MGF) or, equivalently, the cumulant generating function (CGF), of a random variable. Let $X$ be an absolutely continuous random variable with density $f(x)$ with respect to the Lebesgue measure, moment generating function $M(t)$ and CGF $K(t) = \log M(t)$. Then the first-order saddlepoint density approximation to $f(x)$ is given by

$$\hat{f}(x) = \{2\pi K''(\hat{t})\}^{-1/2} \exp\{K(\hat{t}) - \hat{t}x\},$$

where $t = \hat{t}$ is the (unique) solution to the saddlepoint equation $K'(t) = x$, and primes denote derivatives. The Lugannani and Rice (1980) saddlepoint approximation to the cumulative distribution function (CDF) $F(y) = P(X \leq y)$ is obtained by taking $\theta = 0$, $F = F_0$ and $K = K_0$ in (19).

More recent developments include saddlepoint approximations for nonlinear statistics. See Daniels and Young (1991), DiCiccio and Martin (1991) and Jing and Robinson (1994) for further details, and see Jensen (1995) for a rigorous account of the underlying mathematical theory of saddlepoint methods. An extensive discussion of saddlepoint methods and their application will appear in Butler (2004).

Unlike much of this previous work, the current paper uses saddlepoint methods to approximate MGFs of truncated distributions with the view that these approximate MGFs may be used for further saddlepoint inversion. The work is therefore similar in spirit to Fraser, Reid and Wong (1991) and Butler and Wood (2002a).

1.2. *Truncation.* Suppose that $X_i$ denotes a random variable with known MGF $M_i(\theta)$ for $i = 1, \ldots, n$, and that for each $i$ we observe $Y_i = X_i | X_i \in (a_i, b_i)$, that is, $Y_i$ is $X_i$ conditioned to lie in the interval $(a_i, b_i)$. In this paper we are concerned with the following question: *is there a convenient and accurate way to approximate the CGF of $Y_i$ using only $K_i(\theta)$, the CGF of the untruncated variable $X_i$?*

If we are just interested in a single random variable, $Y_1$ say, then the question is probably not of much interest because the density and CDF of $Y_1$ can be expressed simply in terms of the density and CDF of $X_1$, with the latter approximated using the saddlepoint approximations indicated above. However, there are situations in which approximations to the CGFs of the $\{Y_i\}$ are potentially very useful. We mention two such examples.

1.2.1. *Computation of Dirichlet probabilities.* We may wish to construct a saddlepoint approximation for the distribution of the sum $\sum_{i=1}^{n} Y_i$. One such application is to the approximation of certain multivariate CDFs arising in sampling theory and extreme value theory as discussed in Butler and Sutton (1998). For these applications, the multivariate CDF is expressed in terms of the density of $\sum_{i=1}^{n} Y_i$, where the underlying MGFs of the $X_i$ are known. Consider, for example, the probability that an arbitrary Dirichlet vector $\mathbf{D} = (D_1, \ldots, D_n) \sim \text{Dirichlet}\{\gamma = (\gamma_1, \ldots, \gamma_n)\}$ lies in a general rectangular region $(\mathbf{a}, \mathbf{b}) = \prod_{i=1}^{n}(a_i, b_i) \subset (0,1)^n$. If the components of $\mathbf{X} = (X_1, \ldots, X_n)$ are independent with $X_i \sim \text{Gamma}(\gamma_i, 1)$, then the Dirichlet is represented in terms of independent Gammas as $\mathbf{D} = \mathbf{X}/S$, where $S = \sum_{i=1}^{n} X_i$. By independence of $S$ and $\mathbf{X}/S$, the distribution of $\mathbf{D}$ is also the conditional distribution of $\mathbf{X}$ given that $S = 1$. These facts and Bayes' theorem lead to

$$(1) \qquad \Pr\{\mathbf{D} \in (\mathbf{a}, \mathbf{b})\} = f_S\{1 | \mathbf{X} \in (\mathbf{a}, \mathbf{b})\} \frac{\prod_{i=1}^{n} \Pr\{X_i \in (a_i, b_i)\}}{f_S(1)}.$$

Here $f_S\{1 | \mathbf{X} \in (\mathbf{a}, \mathbf{b})\}$ is the density of $Z = \sum_{i=1}^{n} Y_i$ at 1, where $Y_i = X_i | X_i \in (a_i, b_i)$, which we approximate using the saddlepoint density. The other terms are standard computations: $\Pr\{X_i \in (a_i, b_i)\}$ is a gamma probability and $f_S(1)$ is the $\text{Gamma}(\sum_i \gamma_i, 1)$ density of $S$ at 1.

1.2.2. *Ion channel models with time-interval omission.* First, we consider an ion channel model which is represented as a two-state homogeneous semi-Markov process with state space $\{o, c\}$, where state $o$ (state $c$) corresponds to the ion channel being open (closed). Suppose that we observe the process $T_{o,0}, T_{c,1}, T_{o,1}, T_{c,2}, \ldots$, where $T_{o,j}$ is the length of the $j$th sojourn in the open state, and $T_{c,k}$ is the length of the $k$th sojourn in the closed state. We have assumed that the process has started in state $o$; one could equally start in state $c$. Homogeneity and the semi-Markov assumption imply that $\{T_{o,j} : j \geq 1\}$ and $\{T_{c,k} : k \geq 1\}$ are both independent and identically distributed (IID) sequences. Suppose that the MGFs of $T_{o,1}$ and $T_{c,1}$ are, respectively, $\Phi_{oc}(\theta)$ and $\Phi_{co}(\theta)$ and that both are convergent in open neighborhoods of zero.

In ion channel modeling, a phenomenon known as *time interval omission* is commonly built into the model. In effect, this means that only state residences which last for longer than a given time threshold are observed (or detected), and those residences lasting for less than this threshold are not observed (or are undetected); that is, we only observe those sojourns in state $o$ (state $c$) which last at least $\tau_o$ ($\tau_c$); otherwise, it appears to the observer that a jump has not occurred. Time interval omission occurs because of limitations in the sensitivity of the measuring device. Denote the observed sequence by $\widetilde{T}_{o,0}, \widetilde{T}_{c,1}, \widetilde{T}_{o,1}, \ldots$. As a concrete example, suppose that

4       R. W. BUTLER AND A. T. A. WOOD

$T_{c,1} > \tau_c$, $T_{o,1} \leq \tau_o$ and $T_{o,2} > \tau_o$; then $\widetilde{T}_{o,0} = T_{o,0}$, $\widetilde{T}_{c,1} = T_{c,1} + T_{o,1} + T_{c,2}$. The sequences $\{\widetilde{T}_{o,j} : j \geq 1\}$ and $\{\widetilde{T}_{c,k} : k \geq 1\}$ are both IID. Let $\widetilde{\Phi}_{oc}(\theta)$ and $\widetilde{\Phi}_{co}(\theta)$ denote the MGF of a typical member of each sequence. For inferential purposes it is important to express $\widetilde{\Phi}_{oc}$ and $\widetilde{\Phi}_{co}$ in terms of $\Phi_{oc}$ and $\Phi_{co}$.

Define

$$\Phi^D_{oc}(\theta) = E[\exp(\theta T_{o,1}) I(T_{o,1} > \tau_o)], \tag{2}$$

$$\Phi^U_{oc}(\theta) = E[\exp(\theta T_{o,1}) I(T_{o,1} \leq \tau_o)] \tag{3}$$

and

$$\pi_o = P[T_{o,1} > \tau_o] = \Phi^D_{oc}(0),$$

with corresponding definitions for $\Phi^D_{co}$, $\Phi^U_{co}$ and $\pi_c$. Elementary arguments show that $\widetilde{\Phi}_{oc}(\theta)$ can be expressed in terms of a geometric series:

$$\widetilde{\Phi}_{oc}(\theta) = \pi_o^{-1} \Phi^D_{oc} \sum_{n=0}^{\infty} \{\Phi^U_{co}(\theta) \Phi_{oc}(\theta)\}^n \pi_c$$
$$= \pi_o^{-1} \Phi^D_{oc} \{1 - \Phi^U_{co}(\theta) \Phi_{oc}(\theta)\}^{-1} \pi_c.$$

A similar argument shows that

$$\widetilde{\Phi}_{co}(\theta) = \pi_c^{-1} \Phi^D_{co} \{1 - \Phi^U_{oc}(\theta) \Phi_{co}(\theta)\}^{-1} \pi_o.$$

The above discussion shows that time interval omission leads directly to consideration of MGFs of truncated random variables.

More interesting ion channel models have several open states and several closed states, some of which communicate; see Ball, Milne and Yeo (1991). This leads to a more complicated structure, due to aggregation, in which $\Phi_{oc}(\theta)$ and $\Phi_{co}(\theta)$ now represent matrices, each component of which is essentially an MGF which can be expressed as a rational function of the MGFs of the underlying distributions of transition times between individual states. These rational functions are difficult to write down explicitly, but they are straightforward to compute numerically using matrix algebra; see Ball, Milne and Yeo [(1991), Section 3] and also Butler (2000) for analogous calculations in a reliability context. The key point is that, in multistate ion channel models, $\widetilde{\Phi}_{oc}$ and $\widetilde{\Phi}_{co}$ are matrices rather than real numbers, but have similar form to that given above, and are expressed in terms of the MGFs of truncated random variables, as in (2) and (3); see Ball, Milne and Yeo [(1991), Section 4]. Accurate approximation of these distributions can be performed using the methods developed in this paper, but seems to be very difficult otherwise (except in the Markov case).

The present paper was motivated by the ion channel application described above. Further details of this application will be presented in Ball, Butler and Wood (2004).

1.3. *Outline of the paper.* In Section 2 we consider two types of representation for the MGF of a truncated random variable expressed in terms of the MGF of the underlying random variable. One of these representations is obtained by exponential tilting. A second type of representation, which has two versions, is referred to as an exponential convolution representation. In Section 3 we consider saddlepoint approximations to the MGF of the truncated random variable which are motivated by these representations, and indicate their performance in a number of examples. In Section 4, results concerning the tail behavior of the various approximations are given. Proofs of the theorems are given in the Appendix. The research report Butler and Wood (2002b) presents extensions to the lattice and multivariate cases, as well as additional numerical examples.

It turns out that each of the three approximations is extremely accurate in those cases "to which it is suited." Moreover, there is a simple rule of thumb (see Section 3.4) for deciding which approximation to use in a given case. If this rule is followed, numerical and theoretical results indicate that the resulting hybrid approximation will be extremely accurate.

## 2. Representations of truncated MGFs.

2.1. *Preliminaries.* Let $M_0(\theta)$ denote the MGF and $K_0(\theta) = \log M_0(\theta)$ the CGF of a random variable $X$ on $\mathbf{R}$ with density $f_0$ with respect to the Lebesgue measure, and CDF $F_0(x) = P(X \leq x)$. Assume that $M_0(\theta)$ has a convergence strip given by $\theta \in (-\alpha, \beta)$, where $0 < \alpha, \beta \leq \infty$. Let $a < b$ denote real numbers such that $F_0(b) - F_0(a) > 0$.

Let

$$(4) \qquad \mathcal{M}_{(a,b)}(\theta) = \frac{1}{F_0(b) - F_0(a)} \int_a^b e^{\theta x} \, dF_0(x)$$

denote the MGF of $X$ truncated at $a$ and $b$, and conditioned to lie in $(a, b)$. We shall refer to $\mathcal{M}_{(a,b)}(\theta)$ as a truncated MGF which is an abbreviation for the MGF of a truncated random variable, and similar terminology is used for other quantities such as the CGF.

In this paper we discuss how to approximate the truncated CGF $\mathcal{K}_{(a,b)}(\theta) = \log \mathcal{M}_{(a,b)}(\theta)$ and its derivatives in terms of the original CGF $K_0(\theta) = \log M_0(\theta)$ and its derivatives.

2.2. *Tilted representation.* Let $F_\theta(x)$ denote the CDF of the $\theta$-tilted distribution of $X$, that is, $dF_\theta(x) = f_\theta(x) \, dx = e^{\theta x} \, dF_0(x)/M_0(\theta)$, where $f_\theta(x) = e^{\theta x} f_0(x)/M_0(\theta)$ is the density corresponding to the CDF $F_\theta$. Then, for $\theta \in (-\alpha, \beta)$, elementary manipulations show that

$$(5) \qquad \mathcal{M}_{(a,b)}(\theta) = M_0(\theta)[\{F_\theta(b) - F_\theta(a)\}/\{F_0(b) - F_0(a)\}].$$

We shall refer to (5) as the tilted representation of $\mathcal{M}_{(a,b)}(\theta)$.

2.3. *Exponential convolution representations.* We now provide alternative representations of (4) which are collectively valid for all $\theta$ in the convergence interval of $\mathcal{M}_{(a,b)}$. Define

(6) $\quad \Xi_1(\theta, y) = \dfrac{1}{2\pi i} \displaystyle\int_{c_1-i\infty}^{c_1+i\infty} M_0(s) \dfrac{e^{(\theta-s)y}}{\theta - s}\, ds, \qquad -\alpha < c_1 < \min(\beta, \theta)$

and

(7) $\quad \Xi_2(\theta, y) = \dfrac{1}{2\pi i} \displaystyle\int_{c_2-i\infty}^{c_2+i\infty} M_0(s) \dfrac{e^{(\theta-s)y}}{s - \theta}\, ds, \qquad \max(-\alpha, \theta) < c_2 < \beta.$

THEOREM 2.1 (Properties of $\Xi_1$ and $\Xi_2$). *Suppose that $F_0$ is absolutely continuous with density $f_0$, and assume that for some $c \in (-\alpha, \beta)$, there exists a $\nu(c) \in (0, \infty)$ such that*

(8) $\qquad\qquad \displaystyle\int_{t \in \mathbf{R}} |M_0(c+it)|^{1+\nu(c)}\, dt < \infty.$

*Then the following results hold:*

(i) *We have*

(9) $\qquad\qquad \Xi_1(\theta, y) = \displaystyle\int_{-\infty}^{y} e^{\theta x} f_0(x)\, dx, \qquad \theta \in (-\alpha, \infty)$

*and*

(10) $\qquad\qquad \Xi_2(\theta, y) = \displaystyle\int_{y}^{\infty} e^{\theta x} f_0(x)\, dx, \qquad \theta \in (-\infty, \beta).$

*Hence,*

(11) $\qquad\qquad \Xi_1(\theta, y) + \Xi_2(\theta, y) = M_0(\theta), \qquad \theta \in (-\alpha, \beta).$

(ii) *Let $X$ denote a random variable with MGF $M_0(\theta)$ and let $E$ denote an exponential random variable with rate parameter $|\theta|$ which is independent of $X$. When $\theta > 0$,*

(12) $\qquad\qquad \Xi_1(\theta, y) = \theta^{-1} e^{\theta y} f_{X+E}(y);$

*and when $\theta < 0$,*

(13) $\qquad\qquad \Xi_1(\theta, y) = M_0(\theta) - |\theta|^{-1} e^{\theta y} f_{X-E}(y).$

*In the statement of this theorem $f_Z$ denotes the density of a random variable $Z$.*

(iii) *When $\theta > 0$,*

(14) $\qquad\qquad \Xi_2(\theta, y) = M_0(\theta) - \theta^{-1} e^{\theta y} f_{X+E}(y);$

*and when $\theta < 0$,*

(15) $\qquad\qquad \Xi_2(\theta, y) = |\theta|^{-1} e^{\theta y} f_{X-E}(y).$

(iv) *In the respective domains of definition for $\Xi_1$ and $\Xi_2$,*

$$\mathcal{M}_{(-\infty,y)}(\theta) = \Xi_1(\theta, y)/F_0(y) \quad \text{and} \tag{16}$$
$$\mathcal{M}_{(y,\infty)}(\theta) = \Xi_2(\theta, y)/\{1 - F_0(y)\}.$$

(v) *For a general interval $(a,b)$, $M_{(a,b)}(\theta)$ has the alternative representations*

$$(17) \quad \mathcal{M}_{(a,b)}(\theta) = \{\Xi_1(\theta, b) - \Xi_1(\theta, a)\}/\{F_0(b) - F_0(a)\}, \qquad \theta \in (-\alpha, \infty),$$

*and*

$$(18) \quad \mathcal{M}_{(a,b)}(\theta) = \{\Xi_2(\theta, a) - \Xi_2(\theta, b)\}/\{F_0(b) - F_0(a)\}, \qquad \theta \in (-\infty, \beta).$$

We refer to (16)–(18) as exponential convolution representations of the corresponding truncated MGFs.

REMARK 2.1. Condition (8) is a mild smoothness requirement on the underlying density $f_0$. Note that if, for some $c$, (8) holds with $\nu(c) \in (0,1]$, then absolute continuity of $F_0$ follows; see Theorem 11.6.1 in Kawata (1972). However, if we must take $\nu(c) > 1$ for all $c$, then $F_0$ need not be absolutely continuous; see Theorem 13.4.2 in Kawata (1972) for a counterexample.

REMARK 2.2. Although (11) follows immediately from the addition of (9) and (10), it is also interesting to note that (11) is a consequence of Cauchy's theorem; see Butler and Wood [(2002b), Section 2].

**3. Approximations.** We now present approximations to the truncated CGF $\mathcal{K}_{(a,b)}(\theta) = \log \mathcal{M}_{(a,b)}(\theta)$ and its derivatives, distinguishing between the one-sided cases $a = -\infty$ and $b = y < \infty$, and $a = y > -\infty$ and $b = \infty$, and the two-sided case $a > -\infty$ and $b < \infty$.

3.1. *Lugannani and Rice approximation.* Using the tilted representation of the truncated MGF, we obtain

$$\mathcal{K}_{(-\infty,y)}(\theta) = K_0(\theta) + \log\{F_\theta(y)/F_0(y)\}.$$

We may approximate the $\theta$-tilted CDF $F_\theta(y)$ by applying the Lugannani and Rice approximation with the CGF $K_\theta(s) = K_0(\theta + s) - K_0(\theta)$.

If the convergence strip of $K_0(\theta)$ is $\theta \in (-\alpha, \beta)$ with finite $\beta$, then $\mathcal{K}_{(-\infty,y)}(\theta)$ is defined on the larger set $(-\alpha, \infty)$, but it is not clear how to extend this approximation to $\theta \in [\beta, \infty)$. A simple extension is discussed in Butler and Wood [(2002b), Section 5.2], though it turns out that this extension is unsatisfactory.

The Lugannani and Rice approximation to $F_\theta(y)$ is given by

$$\widehat{F}_\theta(y) = \Phi(w_\theta) + \phi(w_\theta)(w_\theta^{-1} - u_\theta^{-1}), \tag{19}$$

where $\Phi$ and $\phi$ are, respectively, the standard normal CDF and density;

$$w_\theta = \operatorname{sgn}(t_y - \theta)[2\{(t_y - \theta)y - K_0(t_y) + K_0(\theta)\}]^{1/2} \tag{20}$$

and $u_\theta = (t_y - \theta)\{K_0''(t_y)\}^{1/2}$, where $\operatorname{sgn}(x) = -1, 0, 1$ depending on whether $x$ is negative, zero or positive; and $t = t_y$ is the unique solution to the saddlepoint equation $K_0'(t) = y$.

The approximation $\widehat{F}_\theta(y)$ is quite simple to use since it is an explicit function of $\theta$ once $t_y$, the saddlepoint for $\theta = 0$, has been determined; thus, the function $\widehat{\mathcal{K}}_{(-\infty,y)}(\theta)$ is available in explicit form once the single saddlepoint solution $t_y$ has been obtained. To see this, note that the saddlepoint for the tilted distribution $\hat{s}_\theta$ solves

$$K_\theta'(\hat{s}_\theta) = K_0'(\hat{s}_\theta + \theta) = y = K_0'(t_y).$$

By uniqueness of the saddlepoint $\hat{s}_\theta + \theta = t_y$, so that only the computation of $t_y$ is required in order to determine $\{\hat{s}_\theta : \theta \in (-\alpha, \beta)\}$. Thus, the CGF approximation

$$\widehat{\mathcal{K}}_{(-\infty,y)}(\theta) = K_0(\theta) + \log\{\widehat{F}_\theta(y)/\widehat{F}_0(y)\}, \qquad \theta \in (-\alpha, \beta), \tag{21}$$

is explicit in $\theta$.

The first two derivatives of the approximation are given by

$$\widehat{\mathcal{K}}'_{(-\infty,y)}(\theta) = K_0'(\theta) + \{\widehat{F}_\theta(y)\}^{-1}\partial \widehat{F}_\theta(y)/\partial\theta$$

and

$$\widehat{\mathcal{K}}''_{(-\infty,y)}(\theta) = K_0''(\theta) + \{\widehat{F}_\theta(y)\}^{-1}\partial^2 \widehat{F}_\theta(y)/\partial\theta^2 - [\{\widehat{F}_\theta(y)\}^{-1}\partial \widehat{F}_\theta(y)/\partial\theta]^2,$$

where

$$\partial \widehat{F}_\theta(y)/\partial\theta = \phi(w_\theta)[\{y - K_0'(\theta)\}(w_\theta^{-3} - u_\theta^{-1}) - (t_y - \theta)^{-2}\{K_0''(t_y)\}^{-1/2}],$$

and the second partial derivative $\partial^2 \widehat{F}_\theta(y)/\partial\theta^2$ is most easily obtained by numerical differentiation.

In the case of $\mathcal{K}_{(y,\infty)}(\theta)$, for $\theta \in (-\alpha, \beta)$, we have the approximations

$$\widehat{\mathcal{K}}_{(y,\infty)} = K_0(\theta) + \log[\{1 - \widehat{F}_\theta(y)\}/\{1 - \widehat{F}_0(y)\}], \tag{22}$$

$$\widehat{\mathcal{K}}'_{(y,\infty)}(\theta) = K_0'(\theta) - \{1 - \widehat{F}_\theta(y)\}^{-1}\partial \widehat{F}_\theta(y)/\partial\theta$$

and $\widehat{\mathcal{K}}''_{(y,\infty)}(\theta) = K_0''(\theta) - T(\theta, y)$, where

$$T(\theta, y) = \{1 - \widehat{F}_\theta(y)\}^{-1}\partial^2 \widehat{F}_\theta(y)/\partial\theta^2 + [\{1 - \widehat{F}_\theta(y)\}^{-1}\partial \widehat{F}_\theta(y)/\partial\theta]^2,$$

with the partial derivatives of $\widehat{F}_\theta(y)$ the same as before.

For general $a < b$, we may approximate $\mathcal{K}_{(a,b)}(\theta) = \log \mathcal{M}_{(a,b)}(\theta)$ by

$$\widehat{\mathcal{K}}_{(a,b)}(\theta) = K_0(\theta) + \log[\{\widehat{F}_\theta(b) - \widehat{F}_\theta(a)\}/\{\widehat{F}_0(b) - \widehat{F}_0(a)\}],$$
(23)
$$\theta \in (-\alpha, \beta).$$

This is an explicit expression in $\theta$ once two saddlepoints have been determined by solving $K_0'(t_a) = a$ and $K_0'(t_b) = b$.

The Lugannani and Rice approximation is exact when applied to the CDF of an arbitrary normal distribution. Therefore, (23) is exact in this case.

3.2. *The exponential convolution approximations.* These approximations are obtained by applying saddlepoint approximations to the integrals defining $\Xi_1(\theta, y)$ and $\Xi_2(\theta, y)$. Denote these saddlepoint approximations by $\widehat{\Xi}_1(\theta, y)$ and $\widehat{\Xi}_2(\theta, y)$. Then in this approach the CGFs $\mathcal{K}_{(-\infty,y)}(\theta)$ and $\mathcal{K}_{(y,\infty)}(\theta)$ are approximated by

(24) $$\widetilde{\mathcal{K}}_{(-\infty,y)}(\theta) = \log\{\widehat{\Xi}_1(\theta, y)/\widehat{\Xi}_1(0, y)\}, \qquad \theta > -\alpha,$$

and

(25) $$\widetilde{\mathcal{K}}_{(y,\infty)}(\theta) = \log\{\widehat{\Xi}_2(\theta, y)/\widehat{\Xi}_2(0, y)\}, \qquad \theta < \beta.$$

To reduce the number of formulae in this section, we shall use the subscripts 1 and 2 to indicate the intervals $(-\infty, y)$ and $(y, \infty)$, respectively.

The saddlepoint approximations to $\Xi_j(\theta, y)$ $(j = 1, 2)$ are given by

(26) $$\widehat{\Xi}_j(\theta, y) = [2\pi\{1 + (\theta - s_{j,\theta})^2 K_0''(s_{j,\theta})\}]^{-1/2} \exp\{K_0(s_{j,\theta}) - (s_{j,\theta} - \theta)y\},$$

where $s_{j,\theta}$ is the unique solution to

(27) $$K_0'(s) + \{\theta - s\}^{-1} = y$$

in $(-\alpha, \beta)$ which satisfies $s_{1,\theta} < \theta$ $(j = 1)$ and $s_{2,\theta} > \theta$ $(j = 2)$.

After some simplifications, we obtain

(28) $$\widetilde{\mathcal{K}}_j(\theta) = \theta y + D_j(\theta) - D_j(0) + K_0(s_{j,\theta}) - K_0(s_{j,0}) - (s_{j,\theta} - s_{j,0})y,$$

where, using implicit differentiation, we have

$$D_j(\theta) \equiv \tfrac{1}{2}\log(\partial s_{j,\theta}/\partial \theta) = -\tfrac{1}{2}\log\{1 + (\theta - s_{j,\theta})^2 K_0''(s_{j,\theta})\}.$$

Note that the approximations are calibrated so that $\widetilde{\mathcal{K}}_j(0) = \mathcal{K}_j(0) = 0$, $j = 1, 2$.

The first derivative of $\widetilde{\mathcal{K}}_j(\theta)$ $(j = 1, 2)$ is given by

$$\widetilde{\mathcal{K}}_j'(\theta) = y + D_j'(\theta) - \{y - K_0'(s_{j,\theta})\}\partial s_{j,\theta}/\partial\theta,$$

where

$$D'_j(\theta) = \tfrac{1}{2}(\partial^2 s_{j,\theta}/\partial\theta^2)/(\partial s_{j,\theta}/\partial\theta), \qquad \partial s_{j,\theta}/\partial\theta = \{1 + (\theta - s_{j,\theta})^2 K_0''(s_{j,\theta})\}^{-1}$$

and

$$\frac{\partial^2 s_{j,\theta}}{\partial\theta^2} = -\frac{(\theta - s_{j,\theta})^2 [K_0'''(s_{j,\theta}) + 2(\theta - s_{j,\theta})\{K_0''(s_{j,\theta})\}^2]}{\{1 + (\theta - s_{j,\theta})^2 K_0''(s_{j,\theta})\}^3}.$$

The second partial derivative $\partial^2 \widetilde{\mathcal{K}}_j(\theta)/\partial\theta^2$ can be determined using numerical differentiation.

In some examples we considered the second-order saddlepoint approximation to $\Xi_j$ given by $\breve{\Xi}_j(\theta, y) = \widehat{\Xi}_j(\theta, y) R(\theta, y)$, where $R(\theta, y)$ is the usual second-order term given in this case by

$$R(\theta, y) = 1 + \frac{1}{8} \frac{K_0''''(s_{j,\theta}) + 6(\theta - s_{j,\theta})^{-4}}{\{K_0''(s_{j,\theta}) + (\theta - s_{j,\theta})^{-2}\}^2} - \frac{5}{24} \frac{\{K_0'''(s_{j,\theta}) + 2(\theta - s_{j,\theta})^{-3}\}^2}{\{K_0''(s_{j,\theta}) + (\theta - s_{j,\theta})^{-2}\}^3}.$$

The resulting approximations to $\mathcal{K}_{(a,b)}(\theta)$ based on (17) and (18) are

(29) $\quad \widetilde{\mathcal{K}}_{1,(a,b)}(\theta) = \log[\{\breve{\Xi}_1(\theta, b) - \breve{\Xi}_1(\theta, a)\}/\{\breve{\Xi}_1(0, b) - \breve{\Xi}_1(0, a)\}]$

for $\theta \in (-\alpha, \infty)$, and for $\theta \in (-\infty, \beta)$,

(30) $\quad \widetilde{\mathcal{K}}_{2,(a,b)}(\theta) = \log[\{\breve{\Xi}_2(\theta, a) - \breve{\Xi}_2(\theta, b)\}/\{\breve{\Xi}_2(0, a) - \breve{\Xi}_2(0, b)\}].$

3.3. *Summary of numerical results.* We now discuss several examples which have been chosen to illustrate some general points. A more extensive set of examples is given in Butler and Wood (2002b). As before, the truncation occurs at $-\infty \leq a < b \leq \infty$ and the convergence strip of the underlying CGF $K_0$ is $(-\alpha, \beta)$, where $0 < \alpha, \beta \leq \infty$.

1. In cases where truncation leads to an extension of the convergence strip of the MGF (i.e., if either $a > -\infty$ and $\alpha < \infty$, or $b < \infty$ and $\beta < \infty$, or both) the most obvious way to extend the LR-based approximation of Section 3.1 is described in Butler and Wood [(2002b), Section 5.2]. However, this extended approximation is poor, as can be seen in Figure 1. The discussion in Butler and Wood [(2002b), Section 5.2] indicates that this is a general problem and not specific to this example.

2. Theoretical results (see Theorem 4.2 and Section 5.3) indicate that when an exponential convolution approximation is used it is appropriate to use (24) or (29) in the right tail and (25) or (30) in the left tail. These findings are strongly supported by our numerical examples; see Figure 2 for a typical case.

3. In our numerical examples, we have found that the first-order saddlepoint approximation to $\Xi_j$ works better in the case of one-sided truncation (i.e., if either $a = -\infty$ or $b = \infty$), while the second-order approximation works better in the case of two-sided truncation (i.e., if $a > -\infty$ and $b < \infty$).

4. In those cases where the convergence strip does not need to be extended, the LR-based approximation has generally proved more accurate than the appropriate exponential convolution approximation, though the latter performs respectably. Figures 3 and 4 present a typical example of this finding.

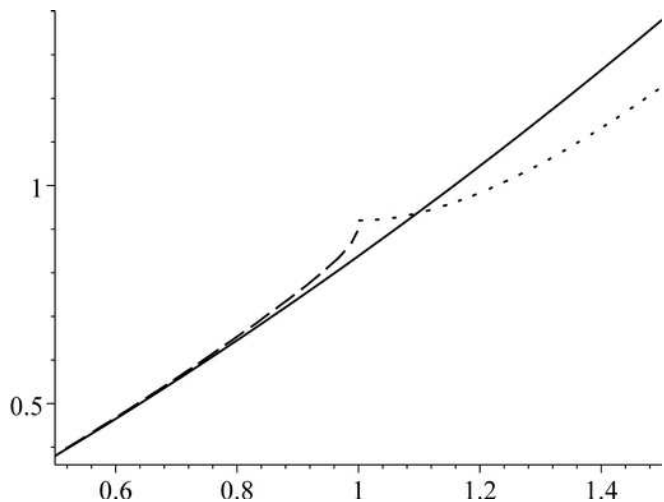

Fig. 1. *Right truncation of an Exponential*(1) *distribution. Plot of* $\mathcal{K}_{(0,2)}(\theta)$ *(solid) and its approximation* $\widehat{\mathcal{K}}_{(0,2)}(\theta)$ *(dashed) for* $\theta \leq 1$ *and its continuation (dotted) for* $\theta \geq 1$.

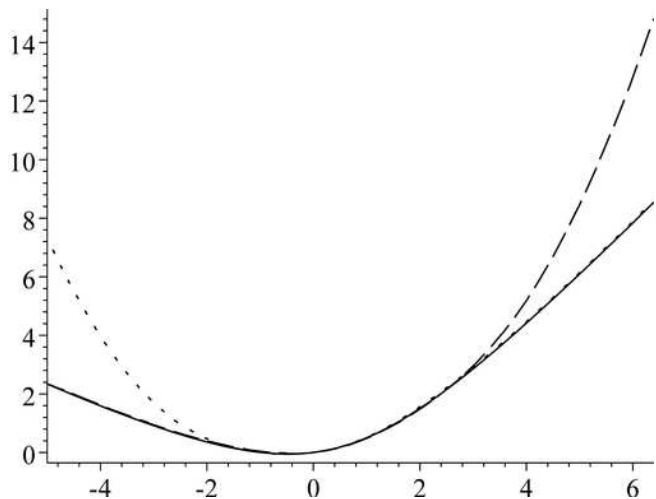

Fig. 2. *Two-sided truncation of Normal*(0,1). *Plot of* $\widetilde{\mathcal{K}}_{1,(-1,2)}(\theta)$ *(dotted)*, $\widetilde{\mathcal{K}}_{2,(-1,2)}(\theta)$ *(dashed)*, *and* $\mathcal{K}_{(-1,2)}(\theta)$ *(solid) for* $\theta \in (-5,6)$.

5. Finally, we return to Example 1.1. The question of interest here is how accurately we can approximate rectangular Dirichlet probabilities using the truncated MGF approximations described above, thereby avoiding the exact computation of the truncated CGF, which is difficult. Table 1 presents results for particular examples. The exponential convolution approximations show consistent accuracy when the saddlepoint is positive; and the Lugannani and Rice approximations are consistently accurate when the saddlepoints are negative. Inaccuracy only arises when either approximation is used in the inappropriate tail.

3.4. *Rule of thumb.* The results of Section 3.3 suggest the following rule of thumb for choosing the approximations, which has worked very well in all the examples we have looked at. In the rule, left and right tail refer to $\theta < 0$ and $\theta \geq 0$, respectively.

*Approximation for right truncation* $(-\infty, y)$. Use the Lugannani and Rice approximation (21) in both tails with one exception. If $\beta < \infty$, so the convergence strip is extended in the right tail, then use (24) in the right tail.

*Approximation for left truncation* $(y, \infty)$. Use the Lugannani and Rice approximation (22) in both tails if there is no extension in the left tail. With extension due to $\alpha > -\infty$, use (25) in the left tail instead.

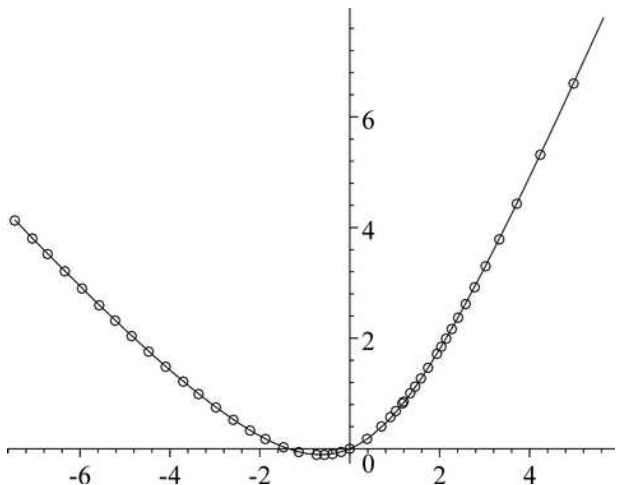

FIG. 3. *Two-sided truncation of the Gumbel distribution. Plot of $\mathcal{K}_{(-1,2)}(\theta)$ (solid) and the "rule of thumb" approximation (centers of circles) that uses the Lugannani and Rice approximation* (23) *for $\theta \leq 0$ and the exponentially convoluted approximation* (29) *for $\theta > 0$.*

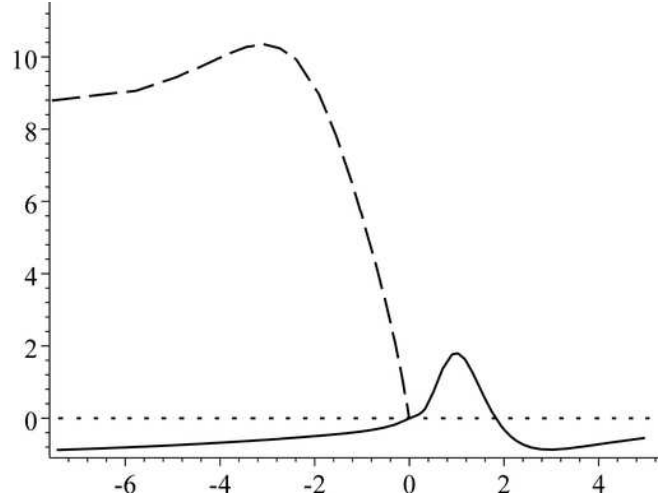

FIG. 4. *Plot of* $100\times$ *absolute error for the "rule of thumb" approximation (solid) in Figure* 3. *If the exponentially convoluted approximation* (30) *replaces* (23) *for* $\theta \leq 0$, *then the error is shown as the dashed line.*

*Approximation for two-sided truncation* $(a, b)$. Use the Lugannani and Rice approximation (23) in those tails in which there are no extensions. Where extensions occur in the left and/or right tails, use (30) and/or (29), respectively.

Since all the approximations to the truncated CGF are calibrated to be zero at $\theta = 0$, it follows that the approximation obtained by following the rule of thumb will be continuous but, in general, not continuously differentiable at $\theta = 0$. However, we have not found the lack of smoothness at $\theta = 0$ to be an issue in practice.

**4. Theoretical accuracy in the tails.** We now investigate the behavior of the approximations to $\mathcal{M}_{(a,b)}(\theta)$ and $\mathcal{K}^{(r)}_{(a,b)}(\theta)$, $r = 1, 2$, as $|\theta| \to \infty$. We make the following assumptions throughout this section:

(A1) The exponential family $\{F_\theta : \theta \in (-\alpha, \beta)\}$ is *steep*, that is, $|K'_0(\theta)| \to \infty$ as $\theta \downarrow -\alpha$ and as $\theta \uparrow \beta$.
(A2) The density $f_0$ has one-sided limits at the truncation points $a$ and $b$, that is, the limits $\lim_{\epsilon \downarrow 0} f_0(a+\epsilon) = f_0(a+)$ and $\lim_{\epsilon \downarrow 0} f_0(b-\epsilon) = f_0(b-)$ both exist.

Note that under (A1) and (A2) and regardless of the value of $a \geq -\infty$ we have, as $\theta \to \infty$,

$$\mathcal{M}_{(a,b)}(\theta) \sim \begin{cases} \theta^{-1} e^{\theta b} f_0(b-)/[F_0(b) - F_0(a)], & \text{if } b < \infty, \\ M_0(\theta)/[1 - F_0(a)], & \text{if } b = \infty, \end{cases}$$

TABLE 1
*Dirichlet probability computations; see Example* 1.1. *For the various values of* $n, \gamma$ *and* $(\mathbf{a}, \mathbf{b})$, *the "Exact" probability as listed was computed using symbolic computation in Maple* V. *The mean of* $Z$ *is listed in the cell "Mean" and its value relative to value* 1 *determines whether the listed saddlepoints for methods* $\widetilde{K}_Z$ *and* $\widehat{K}_Z$ *are negative or positive. Category "SA,* $\widetilde{K}_Z$ *using* $\check{\Xi}_1$*" approximates the CGF of each* $Y_i$ *by using the appropriate (one-sided or two-sided) second-order exponential convolution approximation based on* $\check{\Xi}_1$ *given in* (24) *and* (29), *respectively. Upon determination of* $\widetilde{K}_Z$, *the results of its first-order saddlepoint density inversions are listed. The final column "SA,* $\widehat{K}_Z$ *using L&R" shows comparable computations using the LR-based approximation* $\widehat{K}_Z$ *given by* (21) *or* (23), *as appropriate*

| $n$ | $\gamma$ / Mean | $\mathbf{a}$ / $\mathbf{b}$ | Exact | SA, $\widetilde{K}_Z$ using $\check{\Xi}_1$ Saddlept. | SA, $\widehat{K}_Z$ using L&R Saddlept. |
|---|---|---|---|---|---|
| 3 | $(10, 8, 8)$ | $(0)^3$ | 0.9527 | 0.8877 | 0.9756 |
|   | 1.454 | $\frac{11}{19}, \frac{10}{19}, \frac{11}{19}$ |  | $-23.6$ | $-23.5$ |
| 3 | $(10, 8, 8)$ | $(0)^3$ | 0.02400 | 0.02479 | 0.00141 |
|   | 0.9435 | $0.45, (0.3)^2$ |  | 36.3 | 0.9667 |
| 3 | $(1)^3$ | $(0.2)^3$ | 0.04000 | 0.03869 | 0.0117 |
|   | 0.9080 | $(0.4)^3$ |  | 11.50 | 10.26 |
| 3 | $(1)^3$ | $(0)^3$ | 0.04000 | 0.04059 | 0.02831 |
|   | 0.5707 | $(0.4)^3$ |  | 15.59 | 0.9203 |
| 5 | $(1)^5$ | $(0)^5$ | 0.3680 | 0.3540 | 0.2535 |
|   | 0.9268 | $(0.4)^5$ |  | 1.00 | 0.1335 |
| 5 | $(1, 2, \ldots, 5)$ | $(0)^5$ | 0.5526 | 0.5469 | 0.5439 |
|   | 1.389 | $(0.4)^5$ |  | $-10.14$ | $-9.93$ |
| 5 | $(1, 2, \ldots, 5)$ | $(0)^5$ | $0.0^6 2288$ | $0.0^6 2336$ | $0.0^7 5733$ |
|   | 0.7221 | $0.5, 0.4, (0.1)^3$ |  | 12.65 | 0.9140 |
| 5 | $(1, 2, \ldots, 5)$ | $(0.1)^5$ | 0.03220 | 0.03217 | 0.03183 |
|   | 1.125 | $(0.3)^5$ |  | $-10.02$ | $-9.15$ |
| 10 | $(1)^{10}$ | $(1/15)^5$ | $0.0^4 5080$ | $0.0^4 3175$ | $0.0^4 4969$ |
|   | 1.776 | $(0.3)^5$ |  | $-29.56$ | $-28.56$ |

and regardless of the value of $b \leq \infty$ we have, as $\theta \to -\infty$,

$$\mathcal{M}_{(a,b)}(\theta) \sim \begin{cases} \theta^{-1} e^{\theta a} f_0(a+) / [F_0(b) - F_0(a)], & \text{if } a > -\infty, \\ M_0(\theta) / F_0(b), & \text{if } a = -\infty. \end{cases}$$

4.1. *Accuracy of the Lugannani and Rice approximation.* We first consider the accuracy in the tails of the Lugannani and Rice (LR) approximation $\widehat{\mathcal{M}}_{(a,b)}$ and its logarithmic derivatives $\widehat{\mathcal{K}}'_{(a,b)}$ and $\widehat{\mathcal{K}}''_{(a,b)}$. Theorem 4.1 below is proved in the Appendix.

REMARK 4.1. Comparison of the results in Theorem 4.1 with the limiting results for $\mathcal{M}_{(a,b)}(\theta)$ above shows that the relative error stays bounded in all cases. With $\widehat{\mathcal{K}}'_{(a,b)}(\theta)$ and $\widehat{\mathcal{K}}''_{(a,b)}(\theta)$, the errors actually go to zero as $|\theta| \to \infty$.

THEOREM 4.1. *Consider the LR approximation $\widehat{\mathcal{M}}_{(a,b)}(\theta)$ specified in Section 3.1. Assume that (A1) and (A2) both hold. Suppose also that (i) $\alpha = \infty$ in all statements concerning the left tail and $\beta = \infty$ in all results concerning the right tail; and (ii) as $|\theta| \to \infty$, $u_\theta/w_\theta^3 \to 0$, where $w_\theta$ and $u_\theta$ are given in (20) and below (20), respectively, with $y = a$ or $b$ as appropriate.*

(i) *As $\theta \to \infty$,*

$$\widehat{\mathcal{M}}_{(a,b)}(\theta) \sim \begin{cases} \theta^{-1} e^{\theta b} \hat{f}_0(b-)/[\widehat{F}_0(b) - \widehat{F}_0(a)], & \text{if } b < \infty, \\ M_0(\theta)/[1 - \widehat{F}_0(a)], & \text{if } b = \infty, \end{cases}$$

*and as $\theta \to -\infty$,*

$$\widehat{\mathcal{M}}_{(a,b)}(\theta) \sim \begin{cases} \theta^{-1} e^{\theta a} \hat{f}_0(a+)/[\widehat{F}_0(b) - \widehat{F}_0(a)], & \text{if } a > -\infty, \\ M_0(\theta)/\widehat{F}_0(b), & \text{if } a = -\infty, \end{cases}$$

*where $\hat{f}_0$ is the saddlepoint density approximation to $f_0$ and $\widehat{F}_0$ is the Lugannani and Rice approximation to the CDF $F_0$.*

(ii) *As $\theta \to \infty$,*

$$\widehat{\mathcal{K}}'_{(a,b)}(\theta) = \begin{cases} b - \theta^{-1} + o(\theta^{-1}), & \text{if } b < \infty, \\ K'_0(\theta)(1 + o(1)), & \text{if } b = \infty, \end{cases}$$

*and as $\theta \to -\infty$,*

$$\widehat{\mathcal{K}}'_{(a,b)}(\theta) = \begin{cases} a - \theta^{-1} + o(\theta^{-1}), & \text{if } a > -\infty, \\ K'_0(\theta)(1 + o(1)), & \text{if } a = -\infty. \end{cases}$$

(iii) *As $\theta \to \infty$,*

$$\widehat{\mathcal{K}}''_{(a,b)}(\theta) \sim \begin{cases} \theta^{-2}, & \text{if } b < \infty, \\ K''_0(\theta), & \text{if } b = \infty, \end{cases}$$

*and as $\theta \to -\infty$,*

$$\widehat{\mathcal{K}}''_{(a,b)}(\theta) \sim \begin{cases} \theta^{-2}, & \text{if } a > -\infty, \\ K''_0(\theta), & \text{if } a = -\infty. \end{cases}$$

4.2. *Accuracy of the exponential convolution approximation.* For $j = 1, 2$, let $\widehat{\Xi}_j(\theta, y)$ and $\widetilde{\mathcal{K}}_j^{(r)}(\theta)$, $r = 0, 1, 2$, be as in Section 3.2 and define $\widetilde{\mathcal{M}}_j(\theta) = \widehat{\Xi}_j(\theta, y)/\widehat{\Xi}_j(0, y)$. Also, for $-\infty < a < b < \infty$, define

$$\widetilde{\mathcal{M}}_{1,(a,b)}(\theta) = \{\widehat{\Xi}_1(\theta, b) - \widehat{\Xi}_1(\theta, a)\}/\{\widehat{\Xi}_1(0, b) - \widehat{\Xi}_1(0, a)\},$$

$\widetilde{\mathcal{K}}_{1,(a,b)}(\theta) = \log \widetilde{\mathcal{M}}_{1,(a,b)}(\theta)$, with corresponding definitions for $\widetilde{\mathcal{M}}_{2,(a,b)}(\theta)$ and $\widetilde{\mathcal{K}}_{2,(a,b)}(\theta)$.

REMARK 4.2. Comparison of the results in Theorem 4.2 with the limiting results for $\mathcal{M}_{(a,b)}(\theta)$ shows that the relative error stays bounded in all cases. With $\widetilde{\mathcal{K}}'_j(\theta)$ and $\widetilde{\mathcal{K}}''_j(\theta)$, the errors actually go to zero as $\theta \to \pm\infty$ in the cases covered by the theorem.

THEOREM 4.2. *Assume that* (A1) *and* (A2) *both hold*.

(i) *As* $\theta \to \infty$,
$$\widetilde{\mathcal{M}}_1(\theta) \sim \theta^{-1} e^{\theta y} \widehat{f}_0(y-)/\widehat{\Xi}_1(0,y), \qquad \widetilde{\mathcal{K}}'_1(\theta) = y - \theta^{-1} + o(\theta^{-1})$$
*and* $\widetilde{\mathcal{K}}'_1(\theta) \sim \theta^{-2}$.

(ii) *The limiting behavior of* $\widetilde{\mathcal{M}}_2$, $\widetilde{\mathcal{K}}'_2$ *and* $\widetilde{\mathcal{K}}''_2$ *in the lower tail is the same as that of* $\widetilde{\mathcal{M}}_1(\theta)$, $\widetilde{\mathcal{K}}'_1(\theta)$ *and* $\widetilde{\mathcal{K}}''_1(\theta)$ *in the upper tail, as given in part* (i).

(iii) *If* $-\infty < a < b < \infty$, *then as* $\theta \to \infty$,
$$\widetilde{\mathcal{M}}_{1,(a,b)}(\theta) \sim \theta^{-1} e^{\theta b} \widehat{f}_0(b-)/\{\widehat{\Xi}_1(0,b) - \widehat{\Xi}_1(0,a)\},$$
$\widetilde{\mathcal{K}}'_{1,(a,b)}(\theta) = b - \theta^{-1} + o(\theta^{-1})$ *and* $\widetilde{\mathcal{K}}''_{1,(a,b)}(\theta) \sim \theta^{-2}$.

(iv) *If* $\theta \to -\infty$, *then* $\widetilde{\mathcal{M}}_{2,(a,b)}(\theta)$ *and the derivatives of* $\widetilde{\mathcal{K}}_{2,(a,b)}$ *obey results corresponding to those in part* (iii), *but with* $a$ *replacing* $b$.

PROOF. In part (i), the key point to note is that $s_{1,\theta} \to t_y$ as $\theta \to \infty$, and then the proof follows easily. The proof is essentially the same in the other cases. □

4.3. *Behavior in the other tail.* In Theorem 4.2 we described the limiting behavior of $\widetilde{\mathcal{M}}_1(\theta)$ and its logarithmic derivatives as $\theta \to \infty$, and the behavior of $\widetilde{\mathcal{M}}_2(\theta)$ and its logarithmic derivatives as $\theta \to -\infty$. In this section we indicate, without proof, what happens to $\widetilde{\mathcal{M}}_1(\theta)$ and its derivatives when $\theta \to -\infty$. The results for $\widetilde{\mathcal{M}}_2(\theta)$ are similar and are therefore omitted.

If $\lim_{s \downarrow -\alpha} K''_0(s)/[K'_0(s)]^2 \to 0$, then

(31) $$\widetilde{\mathcal{M}}_1(\theta) \sim M_0(\theta)(e/\sqrt{2\pi})/\widehat{\Xi}_1(0,y) \qquad \text{as } \theta \to -\infty.$$

Under the stronger conditions
$$\lim_{s \downarrow -\alpha} K''_0(s)/K'_0(s) \to 0 \quad \text{and} \quad \lim_{s \downarrow -\alpha} K^{(4)}_0(s)/[K'_0(s)]^3 \to 0,$$
we have

(32) $$\widetilde{\mathcal{K}}'_1(\theta) \sim K'_0(\theta);$$

and still stronger conditions are needed to ensure that

(33) $$\widetilde{\mathcal{K}}''_1(\theta) \sim K''_0(\theta).$$

A sufficient condition for (31)–(33) to hold is the following:

(34) for each $j \geq 2$ $\quad \lim_{s \downarrow -\alpha} K_0^{(j)}(s)$ stays bounded.

Note that condition (34) holds for the normal distribution, gamma distribution (in the left tail) and any other distribution which has bounded support on the left. However, in the case of $-X$, where $X$ has a gamma or inverse Gaussian distribution, or if $X$ has a logistic distribution, then $\widetilde{\mathcal{K}}_1'(\theta)$ and $\widetilde{\mathcal{K}}_1''(\theta)$ do not stay bounded as $\theta \downarrow -\alpha$, and (31)–(33) fail to hold.

## APPENDIX

PROOF OF THEOREM 2.1. Using the convolution formula for densities [see, e.g., Theorem 6.1.2 in Chung (1974), for a precise statement], we have, for $\theta > 0$,

$$\int_{-\infty}^{y} e^{\theta x} f_0(x)\, dx = \frac{e^{\theta y}}{\theta} \int_{-\infty}^{\infty} \theta e^{-\theta(y-u)} I(u \leq y) f_0(u)\, du = \frac{e^{\theta y}}{\theta} f_{X+E}(y),$$

where $E$ is an exponential random variable with rate parameter $\theta$, which is independent of $X$. Define

$$H_{c,\theta}(t) = [M_0(c+it)/\{1-(c+it)/\theta\}]/[M_0(c)/\{1-c/\theta\}],$$

so that $H_{0,\theta}(t)$ is the characteristic function (CF) of $f_{X+E}(y)$, and $H_{c,\theta}(t)$ is the CF of the $c$-tilted density $f_{X+E}(y)e^{cy}/\{M_0(c)/(1-c/\theta)\}$. Note that if (8) holds for some $c \in (-\alpha, \beta)$, then (8) holds for all $c \in (-\alpha, \beta)$. For a proof of this result, which involves two applications of the Hausdorff–Young inequality, see Lemma 2.3.4 of Jensen (1995). Using Hölder's inequality,

$$\int_{-\infty}^{\infty} |H_{c,\theta}(t)|\, dt \leq \frac{1-c/\theta}{M_0(c)} \left(\int_{-\infty}^{\infty} |M_0(c+it)|^{1+\nu(c)}\, dt\right)^{1/(1+\nu(c))}$$
$$\times \left(\int_{-\infty}^{\infty} \frac{1}{|1-(c+it)/\theta|^{(1+\nu(c))/\nu(c)}}\, dt\right)^{\nu(c)/(1+\nu(c))} < \infty$$

for each $c \in (-\alpha, \min(\beta, \theta))$. Therefore we may apply the Fourier inversion theorem [see, e.g., Chung (1974), page 155, for a precise statement] to $H_{c,\theta}(t)$ to obtain

(35) $\quad \dfrac{1}{2\pi} \int_{-\infty}^{\infty} H_{c,\theta}(t) e^{-ity}\, dt = f_{X+E}(y) e^{cy}/\{M_0(c)/(1-c/\theta)\}.$

After some rearrangement, we find that (35) gives (9) for all $\theta > 0$ and $c \in (-\alpha, \min(\beta, \theta))$. This shows also that $\Xi_1(\theta, y)$ does not depend on the choice of $c_1$ in (6). An analytic continuation argument extends (9) to $\theta \in (-\alpha, 0]$, thus (9) is established for all $\theta > -\alpha$.

Identical reasoning gives (10) and (15), and (11) follows immediately after adding (9) and (10). The statements (16), (17) and (18) follow directly from the definitions. □

PROOF OF THEOREM 4.1. The LR approximation to $\mathcal{M}_{(a,b)}(\theta)$ is given by

$$\widehat{\mathcal{M}}_{(a,b)}(\theta) = M_0(\theta)[\{\widehat{F}_\theta(b) - \widehat{F}_\theta(a)\}/\{\widehat{F}_0(b) - \widehat{F}_0(a)\}].$$

(i) *Case* $\theta \to \infty$, $b = \infty$. Note that $\widehat{F}_\theta(b) = 1$ for all $\theta$ and $\widehat{F}_\theta(a) \to 0$ as $\theta \to \infty$, so $\widehat{\mathcal{M}}_{(a,b)}(\theta) \sim M_0(\theta)/[1 - \widehat{F}_0(a)]$ as required.

*Case* $\theta \to \infty$, $b < \infty$. Here $\widehat{F}_\theta(a)/\widehat{F}_\theta(b) \to 0$, so

$$\widehat{\mathcal{M}}_{(a,b)}(\theta) \sim M_0(\theta)\widehat{F}_\theta(b)/[\widehat{F}_0(b) - \widehat{F}_0(a)].$$

By assumption $u_\theta/w_\theta^3 \to 0$ as $\theta \to \infty$. Moreover, elementary calculations show that as $w_\theta \to -\infty$, $\Phi(w_\theta) \sim -\phi(w_\theta)[w_\theta^{-1} + w_\theta^{-3}]$, and it then follows easily that

$$\widehat{F}_\theta(b) \sim -\phi(w_\theta)/u_\theta \sim \theta^{-1} e^{b\theta} \hat{f}_0(b-)/M_0(\theta) \qquad \text{as } \theta \to \infty,$$

where $\hat{f}_0(b) = (2\pi)^{-1/2}|K_0''(t_b)|^{-1/2}\exp\{K_0(t_b) - t_b b\}$ is the saddlepoint approximation to $f_0(b)$. The proofs for $\theta \to -\infty$ with $a \geq -\infty$ are similar.

(ii) We have

$$\widehat{\mathcal{K}}'_{(a,b)}(\theta) = K_0'(\theta) + \{\widehat{F}_\theta(b) - \widehat{F}_\theta(a)\}^{-1}[\partial \widehat{F}_\theta(b)/\partial \theta - \partial \widehat{F}_\theta(a)/\partial \theta].$$

*Case* $\theta \to \infty, b = \infty$. Since $\widehat{F}_\theta(b) - \widehat{F}_\theta(a) \to 1$, $\partial \widehat{F}_\theta(b)/\partial \theta = 0$, $\partial \widehat{F}_\theta(a)/\partial \theta \to 0$ and $K_0'(\theta) \to \infty$, the result follows.

*Case* $\theta \to \infty$, $b < \infty$. Here

$$\widehat{F}_\theta(a)/\widehat{F}_\theta(b) \to 0 \quad \text{and} \quad (\partial \widehat{F}_\theta(a)/\partial \theta)/(\partial \widehat{F}_\theta(b)/\partial \theta) \to 0.$$

Therefore,

$$\begin{aligned}\widehat{\mathcal{K}}'_{(a,b)}(\theta) &= K_0'(\theta) + \{\widehat{F}_\theta(b)\}^{-1} \partial \widehat{F}_\theta(b)/\partial \theta + o(\theta^{-1}) \\ &= K_0'(\theta) + b - K_0'(\theta) + \theta^{-1} + o(\theta^{-1}) \\ &= b + \theta^{-1} + o(\theta^{-1}),\end{aligned}$$

as required. The cases $\theta \to -\infty$ with $a = -\infty$ and $a > -\infty$ are proved in similar fashion.

(iii) The results here follow from similar but more extensive calculations. □

Department of Statistics  
Colorado State University  
Fort Collins, Colorado 80523  
USA  
e-mail: walrus@stat.colostate.edu

School of Mathematical Sciences  
University of Nottingham  
Nottingham NG7 2RD  
United Kingdom  
e-mail: atw@maths.nott.ac.uk